\documentclass[letterpaper, 10 pt, conference]{ieeeconf}

\IEEEoverridecommandlockouts

\usepackage{graphics}
\usepackage{epsfig}
\usepackage{times}
\usepackage{amsmath}
\usepackage{amssymb}
\usepackage{multirow}
\usepackage{array}
\usepackage{amsfonts}
\usepackage{subcaption}

\usepackage{amsthm}
\usepackage{todonotes}
\usepackage{comment}
\usepackage{cite}
\usepackage{booktabs}
 \usepackage[colorlinks=true]{hyperref}
\hypersetup{
    colorlinks=true,
    linkcolor=blue,
    filecolor=blue,      
    urlcolor=blue,
    pdftitle={Overleaf Example},
    pdfpagemode=FullScreen,
    }

\theoremstyle{definition}
\newtheorem{assumption}{Assumption}
\newtheorem{definition}{Definition}
\theoremstyle{plain}

\newtheorem{remark}{Remark}
\newtheorem{problem}{Problem}
\theoremstyle{definition}

\title{\LARGE \bf
Platoon Formation in a Mixed Traffic Environment: A Model-Agnostic Optimal Control Approach}

\author{A M Ishtiaque Mahbub, {\itshape{Student Member, IEEE}}, Andreas A. Malikopoulos, {\itshape{Senior Member, IEEE}}
\thanks{This research was supported by ARPAE's NEXTCAR program under the award number DE-AR0000796.}
\thanks{The authors are with the Department of Mechanical Engineering, University of Delaware, Newark, DE 19716 USA (emails: \tt\small{mahbub@udel.edu};\tt\small{andreas@udel.edu}.)}}

\begin{document}

\maketitle
\thispagestyle{empty}
\pagenumbering{arabic}

\begin{abstract}
Coordination of connected and automated vehicles (CAVs) in a mixed traffic environment poses significant challenges due to the presence of human-driven vehicles (HDVs) with stochastic dynamics and driving behavior. In earlier work, we addressed the problem of platoon formation of HDVs led by a CAV using a model-dependent controller. In this paper, we develop a comprehensive model-agnostic, multi-objective optimal controller which ensures platoon formation by directly controlling the leading CAV without having explicit knowledge of the trailing HDV dynamics. We provide a detailed exposition of the control framework that uses instantaneous motion information from multiple successive HDVs to enforce safety while achieving the optimization objectives. To demonstrate the efficacy of the proposed control framework, we evaluate its performance using numerical simulation and provide associated sensitivity and robustness analysis.

\end{abstract}

\indent


\section{Introduction}
The implementation of an emerging transportation system with connected automated vehicles (CAVs) enables a novel computational framework to better monitor the transportation network conditions and make optimal operating decisions to improve safety and reduce pollution, energy consumption, and travel delays \cite{Guanetti2018}.
Recent efforts have reported several optimal control approaches for coordination of CAVs at different traffic scenarios such as on-ramp merging roadways \cite{Ntousakis:2016aa}, roundabouts \cite{Malikopoulos2018a, bakibillah2019optimal}, speed reduction zones \cite{Malikopoulos2018c}, signal-free intersections \cite{Colombo2014,Kim2014,Au2015,Mahbub2019ACC,Malikopoulos2020}, and traffic corridors \cite{Lee2013,mahbub2020decentralized,Beaver2020DemonstrationCity}. These approaches have focused on 100\% CAV penetration rates without considering human-driven vehicles (HDVs). However, the existence of having a transportation network with a 100\% CAVs is not expected before 2060 \cite{alessandrini2015automatedmixed2060}. Therefore, the need for a mathematically robust and tractable control framework considering a \emph{mixed traffic environment} consisting of both CAVs and HDVs is essential. In reality, HDVs pose significant modeling and control challenges to the CAVs due to the stochastic nature of the human-driving behavior, often emulated by the car-following models, see \cite{Gipps1981, wiedemann1974, treiber2013traffic}. Some approaches reported in the literature \cite{ Malikopoulos2018a,wan2016optimalmixed} have included the car-following models for coordinating CAVs in a mixed environment, while others have been based on reinforcement learning \cite{kreidieh2018dissipating, wu2017framework}. 

In this paper, our research hypothesis is that, since we cannot control the HDVs directly, we can control the CAVs in a way to force the trailing HDVs to form platoons, and thus indirectly control the HDVs. In this context, we focus on the problem of vehicle platoon formation in mixed traffic environment by only controlling the CAVs within the network. Although the problem of platoon formation has been widely studied for $100 \%$ CAV penetration \cite{zheng2017platooning, dunbar2006distributed,  orosz2016connected}, only limited efforts have been reported in the literature for mixed traffic environment. Some of these approaches have adopted adaptive cruise control for the CAVs \cite{jin2017fuel_mixplatoon, hajdu2019robust, dollar2021mpc} to maintain platoon stability. 

In this paper, we extend our previous work \cite{mahbub2021_platoonMixed} by introducing a multi-objective optimal control framework for each CAV within the network subject to its state and control constraints. The optimization objectives of the CAV are (a) to form a platoon with the trailing HDVs, and (b) to improve its fuel economy while achieving (a). Our proposed control framework is \emph{model-agnostic}, i.e., it does not require the explicit knowledge of the HDVs' car-following model, and employs a receding horizon controller that uses a \emph{multi-successor communication topology}, i.e., reception of instantaneous motion information from multiple trailing HDVs, to enforce safety while deriving and implementing the optimal control input of the CAV. 
To the best of our knowledge, such approach has not yet been reported in the literature to date.


The remainder of the paper proceeds as follows. In Section \ref{sec:pf}, we provide the modeling framework of the platoon formation problem. In Section \ref{sec:mpc}, we develop a model-agnostic constrained multi-objective optimal control framework for the CAV for platoon formation. In Section \ref{sec:sim}, we evaluate the performance of the proposed control framework using numerical simulation and validate its effectiveness. Finally, we draw concluding remarks and discuss potential directions for future research in Section \ref{sec:conc}.

\section{Problem Formulation}\label{sec:pf}
\begin{figure}[b]
    \centering
    \includegraphics[scale=0.35]{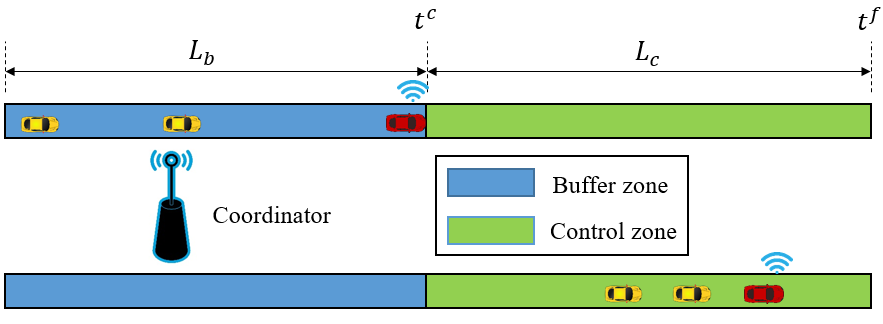}
    \caption{A CAV (red) traveling with two trailing HDVs (yellow), where the HDVs' state information is estimated (top scenario) by the coordinator within the buffer zone, and the platoon is formed (bottom scenario) by controlling the CAV at the control zone.}
    \label{fig:platoon_zone}
\end{figure}

We consider a CAV followed by one or multiple HDVs traveling in a single-lane roadway of length $L\in \mathbb{R}^+$. We subdivide the roadway into a \emph{buffer zone} of length $L_{b}\in \mathbb{R}^+$, where the HDVs' state information is estimated, as shown in Fig. \ref{fig:platoon_zone} (top), and a \emph{control zone} of length $L_{c}\in \mathbb{R}^+$ such that $L=L_{b}+L_{c}$, where the leading CAV is to be controlled to form a platoon with the trailing HDVs, as shown in Fig. \ref{fig:platoon_zone} (bottom). The CAV enters and leaves the control zone at times $t^c, t^f\in \mathbb{R}^+$, respectively.

Let $\mathcal{N}=\{1,\ldots, N\}$, where $N\in \mathbb{N}$ is the total number of vehicles traveling within the buffer zone, be the set of vehicles considered to form a platoon. Here, the leading vehicle indexed by $1$ is the CAV, and the rest of the trailing vehicles in  $\mathcal{N}_{\text{HDV}}:=\mathcal{N}\setminus\{1\}$ are HDVs.  
The objective of CAV $1\in\mathcal{N}$ is to derive and implement the optimal control input (acceleration/deceleration) such that a platoon formation with the trailing HDVs in $\mathcal{N}_{\text{HDV}}$ is completed within the control zone of length $L_c$.

Since the HDVs do not share their local state information with any external agents, we assume the presence of a \emph{coordinator} that gathers the state information of the trailing HDVs traveling within the buffer zone. The coordinator, which can be loop-detectors or comparable sensory devices, in turn, transmits the HDV state information to the CAV at each time instance $t\in[t^c, t^f]$ using standard vehicle-to-infrastructure communication protocol.
\begin{figure}
    \centering
    \includegraphics[scale=0.40]{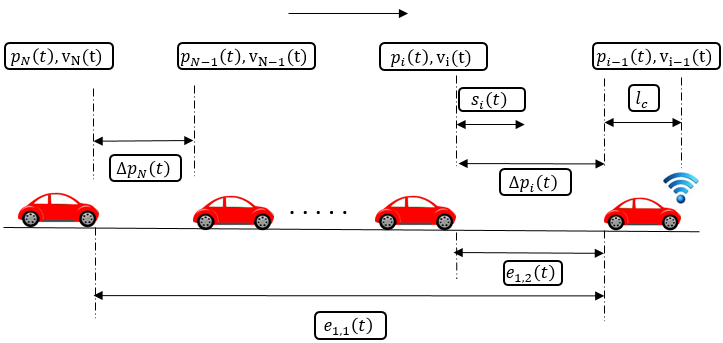}
    \caption{Predecessor-follower coupled car-following dynamic.}
    \label{fig:problem_formulation}
\end{figure}
%
We consider a standard double-integrator model to represent the longitudinal dynamics of each vehicle $i\in\mathcal{N}$ within the network at time $t\in[t^c, t^f]$ as
\small
\begin{align}\label{eq:dynamics_pv}
    \dot{p}_i(t) = v_i(t), ~ 
    \dot{v}_i(t) = u_i(t),
\end{align}
\normalsize
where $p_i(t)\in \mathcal{P}_i$, $v_i(t)\in \mathcal{V}_i$ and $u_i(t)\in\mathcal{U}_i$ are the position of the front bumper, speed and control input (acceleration command) of each vehicle $i\in \mathcal{N}$, respectively.

The speed $v_i(t)$ and control input $u_i(t)$ of each vehicle $i\in \mathcal{N}$ at time $t\in[t^c, t^f]$ are subject to the following constraints,
\small
\begin{align}\label{eq:state_control_constraints}
    0\le v_{\min} \le v_i(t) \le v_{\max},~
    u_{\min} \le u_i(t) \le u_{\max},
\end{align}
\normalsize
where $v_{\min}$ and $v_{\max}$ are the minimum and maximum allowable speed of the considered roadway, respectively, and $u_{\min}$ and $u_{\max}$ are the minimum and maximum acceleration of each vehicle $i\in \mathcal{N}$, respectively.
\begin{definition}
The dynamic following spacing $s_i(t)$ between vehicle $i \text{ and }(i-1)\in\mathcal{N}$ is $ s_i(t)= \rho v_i(t)+ s_0$,
where $\rho\in \mathbb{R}^+$ denotes a desired safety time headway that each HDV $i\in\mathcal{N}_{\text{HDV}}$ maintains while following its preceding vehicle $i-1\in\mathcal{N}$, and $s_0\in \mathbb{R}^+$ is the standstill distance denoting the minimum bumper-to-bumper gap at stop.
\end{definition}
\begin{definition}\label{def:2}
The \textit{headway} $\Delta p_i(t)$ (see Fig. \ref{fig:problem_formulation}) and the \emph{approach rate} $\Delta v_i(t)$ of vehicle $i\in\mathcal{N}$ between two consecutive vehicles $i,~(i-1) \in\mathcal{N}$ are $\Delta p_i(t)=p_{i-1}(t)- p_i(t)-l_c$ and
    $\Delta v_i(t) = v_{i-1}(t) - v_i(t)$, respectively,
where $l_c\in\mathbb{R}^+$ is the length of each vehicle $i\in\mathcal{N}$.
\end{definition}
The rear-end collision avoidance constraint is
\small
\begin{align}\label{eq:rearend_constraint}
    \Delta p_i(t) \ge s_i(t), \quad \forall t\in[t^c, t^f].
\end{align}
\normalsize
%

In our modeling framework, we impose the following assumptions.
\begin{assumption}\label{assum:1}
For each vehicle $i\in\mathcal{N}$ at time $t=t^c$, none of the state, control and safety constraints in \eqref{eq:state_control_constraints}-\eqref{eq:rearend_constraint} are active. 
\end{assumption}
\begin{assumption}\label{assum:4}
No error or delay occurs during the communication between the CAV and the coordinator.
\end{assumption}

Assumption \ref{assum:1} ensures that the initial state and control input of each vehicle $i\in\mathcal{N}$ are feasible.
Assumption \ref{assum:4} may be strong, but it is relatively straightforward to relax as long as the noise in the measurements and/or delays is bounded. 

The control input $u_i(t)$ of each vehicle $i\in\mathcal{N}$ in \eqref{eq:dynamics_pv} can take different forms based on the consideration of connectivity and automation.
For CAV $1\in\mathcal{N}$, we derive and implement the control input $u_1(t)$ using the optimal control framework discussed in Section \ref{sec:mpc}.
%
For each HDV $i\in \mathcal{N}_{\text{HDV}}$, however, we consider a car-following model to represent the predecessor-follower coupled dynamics (see Fig. \ref{fig:problem_formulation}), which has the generic structure $ {{u}_i(t) = f (\Delta p_i(t), \Delta v_i(t), v_i(t))}$.
Here, $f(\cdot)$ represents the behavioral function of the car-following model. In this paper, we consider that the HDVs' behavioral function $f$ is unknown to CAV $1$.  
\begin{definition}\label{def:1}
The information set $\mathcal{I}_1(t)$ of CAV $1\in\mathcal{N}$ at time $t\in[t^c,t^f]$ is $\mathcal{I}_1(t) = \{{p}_{1:N}(t), {v}_{1:N}(t)\}$,
where ${p}_{1:N}(t)=[{p}_1(t),\ldots, {p}_N(t)]^T$ and ${v}_{1:N}(t)=[{v}_1(t),\ldots, {v}_N(t)]^T$.
\end{definition}

\begin{definition}\label{def:uniform_flow}
A platoon formation is established at some time $t^p \in (t^c, t^f)$ if for each vehicle $i\in\mathcal{N}$, the headway $\Delta p_i(t)$ converges to an equilibrium headway $\Delta p_{eq}$, and the approach rate $\Delta v_i(t)$ converges to zero, i.e., 
\small
\begin{align}
 &\Delta p_i(t)= \Delta p_{eq},~\Delta p_{eq}\in \mathbb{R}^+,~ &\forall t\ge t^p,\label{eq:uniform_flow_1}\\
 &{\Delta v_i(t)=0},~ &\forall t\ge t^p,\label{eq:uniform_flow_2}
\end{align}
\normalsize
\end{definition}
\begin{remark}\label{rem:rmse_platoonFormation}
In real-world applications, conditions \eqref{eq:uniform_flow_1}-\eqref{eq:uniform_flow_2} might be too restrictive to establish a platoon formation. Therefore, we relax these conditions and introduce the following root-mean-squared error based conditions to establish a platoon formation at some time $t^p \in (t^c, t^f)$,
\small
\begin{subequations}\label{eq:platoon_formation_rmse}
\begin{align}
    &\sqrt{\frac{1}{N-1}\sum_{i=2}^N ( \Delta p_i(t) - \mu_{\Delta p}(t) )^2} \le \epsilon_{\Delta p},&\forall t\ge t^p,\label{eq:platoon_formation_1}\\
    &\sqrt{\frac{1}{N}\sum_{i=1}^N ( v_i(t) - \mu_{v}(t) )^2} \le \epsilon_{v},&\forall t\ge t^p,\label{eq:platoon_formation_2}
\end{align}
\end{subequations}
\normalsize
where $\mu_{\Delta p}(t):=\frac{\sum_{i=2}^N \Delta p_i(t)}{N}$ and $\mu_{v}(t):= \frac{\sum_{i=1}^N v_i(t)}{N}$ are the mean headway and mean speed of $N$ vehicles, respectively, and  $\epsilon_{\Delta p}$ and $\epsilon_{v}$ are the allowable deviation of $\Delta p_i(t)$ and $v_i(t)$ from the equilibrium values $\Delta p_{eq} \text{ and }v_{eq}$, respectively. 
\end{remark}
%
%
Next, we formally state the platoon formation problem in mixed environment as follows.
\begin{problem}\label{prob:1}
    Given the information set $\mathcal{I}_1(t)$ for each time $t\in[t^c, t^f]$, the objective of the CAV $1\in\mathcal{N}$ is to derive its optimal control input $u_1^*(t)$ so that the each vehicle $i\in\mathcal{N}$ achieves a platoon formation (Definition \ref{def:uniform_flow}) within the control zone.
\end{problem}

\begin{remark}\label{rem:optimization_objectives}
In our framework, CAV $1$ derives its optimal control input $u_1^*(t)$ by solving an optimal control problem with the following objectives: (a) formation of platoon with the trailing HDVs (Definition \ref{def:uniform_flow}), and (b) improvement of its fuel economy while achieving (a).
\end{remark}
%

In this paper, we adopt a receding horizon control framework with multi-successor communication topology to address Problem \ref{prob:1}. In what follows, we provide a detailed exposition of the receding horizon control framework that leads to an optimal platoon formation (Remark \ref{rem:optimization_objectives}).

\section{Receding Horizon Control}\label{sec:mpc}

The basic principle of a receding horizon control is that, the optimal control input sequence at current time instance is obtained by solving an optimal control problem online with the prediction horizon $T_p$, and only implementing the first element of the solved optimal control input sequence. Then the horizon moves forward one step, and the above process is repeated until the optimization horizon $T_h$ is reached. 

\begin{remark}\label{rem: T_h}
The exit time $t^f$ of CAV $1$ from the control zone depends on the nature of the optimal control input of CAV $1$, and thus, it is not known a priori. Let $t^e$ be the time that the CAV exits the control zone when cruising with a constant speed inside the control zone. Then,  $t^e=t^c+\frac{L_c}{v_1(t^c)}$. In our previous work \cite{mahbub2021_platoonMixed}, we have shown that a platoon formation with trailing HDVs can be achieved by non-positive control trajectory of the CAV. Consequently, if we aim at forming the platoon by considering the optimization horizon to be $T_h= t^e-t^c$, then we can ensure that the platoon is formed within the control zone.
\end{remark}
For CAV $1$, we aim to achieve the optimization objectives outlined in Remark \ref{rem:optimization_objectives} while enforcing rear-end collision avoidance constraint with its trailing HDV. To this end, the adoption of the CAV dynamics in \eqref{eq:dynamics_pv} is not sufficient; our proposed control framework requires the consideration of an augmented CAV dynamics model.
\subsection{Augmented CAV dynamics}
To capture the additional characteristics of the platoon formation dynamics from the CAV's control point of view, our proposed control framework uses instantaneous motion information from multiple successive HDVs. Hence, we define two additional states as follows.
\begin{definition}\label{def:error states}
 The \emph{head-to-tail gap} of the platoon, $e_{1,1}(t)$ and the \emph{leader-follower gap}, $e_{1,2}(t)$ are $e_{1,1}(t)= p_1(t)-p_N(t)-(N-1)l_c$ and $e_{1,2}(t)= p_1(t)-p_2(t)-l_c$, respectively (see Fig. \ref{fig:problem_formulation}).
\end{definition}
The additional states  $e_{1,1}(t)$ and $e_{1,2}(t)$ enables the augmentation of the CAV dynamics \eqref{eq:dynamics_pv} with the following set of equations,
\vspace{-0.5cm}
\small
\begin{align}
    &\dot{e}_{1,1}(t) = v_1(t)-v_N(t),\label{eq:error_dynamics1}\\
    &\dot{e}_{1,2}(t) = v_1(t)-v_2(t).\label{eq:error_dynamics2}
\end{align}
\normalsize
\begin{remark}\label{rem:error-states}
The consideration of the {head-to-tail gap} $e_{1,1}(t)$  of the platoon enables the formulation of the objective function for the platoon formation problem whereas the {leader-follower gap} $e_{1,2}(t)$ enables the enforcement of rear-end collision avoidance constraint in \eqref{eq:rearend_constraint}, leading to a safe platoon formation.
\end{remark}
\subsection{Discrete Time Formulation}
To enable the application of discrete time receding horizon control, we formulate the optimal control problem in discrete time. Suppose, the optimization horizon $T_h$ is discretized by a sampling time interval $\tau$ leading to discrete time instance $k$. Assuming constant value of control input $u_1(k)$ during each time step $[k, (k+1)]$, we recast the augmented CAV dynamics \eqref{eq:dynamics_pv} and \eqref{eq:error_dynamics1}-\eqref{eq:error_dynamics2} as linear discrete-time state equations
\small
\begin{align}
    &p_1(k+1) = p_1(k) + v_1(k) \tau + \frac{1}{2}u_1(k) \tau^2,\label{eq:discrete_dynamics_1}\\
    &v_1(k+1) = v_1(k) + u_1(k) \tau,\\
    &e_{1,1} (k+1) = e_{1,1}(k) + (v_1(k)-v_N(k)) \tau + \frac{1}{2}u_1\tau^2,\\
    &e_{1,2} (k+1) = e_{1,2}(k) + (v_1(k)-v_2(k)) \tau  + \frac{1}{2}u_1\tau^2.\label{eq:discrete_dynamics_4}
\end{align}
\normalsize
We define the current state vector $x_1(k)$, measured output vector $y_1(k)$ and the measured disturbance vector $w_1(k)$ as
\small
\begin{gather*}
x_1(k):=
 \begin{bmatrix}
p_1(k)\\ 
v_1(k)\\ 
e_{1,1}(k)\\ 
e_{1,2}(k)
\end{bmatrix},~
y_1(k):=
 \begin{bmatrix}
v_1(k)\\ 
e_{1,1}(k)\\ 
e_{1,2}(k)
\end{bmatrix},~
w_1(k) :=
 \begin{bmatrix}
v_N(k)\\
v_2(k)
\end{bmatrix}.
\end{gather*}
\normalsize
The state-space representation of the discrete dynamic in \eqref{eq:discrete_dynamics_1}-\eqref{eq:discrete_dynamics_4} is thus
\small
\begin{align}
    &x_1(k+1) = Ax_1(k)+ B_u u_1(k) + B_w w_1(k),\label{eq:state_space_1}\\
    &y_1(k) = Cx_1(k),\label{eq:state_space_2}
\end{align}
\normalsize
where, the corresponding state matrix $A$, control matrix $B_u$, disturbance matrix $B_w$ and output matrix $C$ can be computed using \eqref{eq:discrete_dynamics_1}-\eqref{eq:state_space_2}.
%
For the remainder of this paper, we drop the subscript $1$ denoting the CAV from the discrete state-space model where it does not introduce ambiguity.
\subsection{Prediction Model}
In order to solve an online optimization within the prediction horizon $T_p$, the receding horizon controller requires a prediction model to take into account the future possible states. In general, the future system states are predicted based on the model \eqref{eq:state_space_1}-\eqref{eq:state_space_2} and the current state information $x(k)$. 
Let us define the predicted state, predicted output, control and disturbance vector given the prediction horizon $T_p$ and control horizon $T_c$ as $\Tilde{X}(k+T_p|k)=[\Tilde{x}(k+1|k),~\Tilde{x}(k+2|k),\ldots,~\Tilde{x}(k+T_p|k) ]^T$, $\Tilde{Y}(k+T_p|k)=[\Tilde{y}(k+1|k),~\Tilde{y}(k+2|k),\ldots,\Tilde{y}(k+T_p|k)]^T$,~ ${U}(k+{T_c})=[{u}(k),~{u}(k+1),\ldots,{u}(k+T_c-1) ]^T$ and ${W}(k+{T_p})=[{w}(k),~{w}(k+1),\ldots,{w}(k+T_p-1) ]^T$, respectively. Here $\Tilde{x}(k+n|k)$, $\Tilde{y}(k+n|k)$, and ${w}(k+n-1)$, $n=1,\ldots, T_p$, denote the predicted state, output and disturbance values within the prediction horizon $T_p$ based on their value at the discrete instance $k$, respectively. 

The predictive state and associated performance vectors of the receding horizon controller can subsequently be represented as
\small
\begin{gather}
    \Tilde{X}(k+T_p|k) = \Tilde{A}x(k) + \Tilde{B_u}U(k+T_c)+ \Tilde{B_d}W(k+T_p),\label{eq:state_space_pred_1}\\
    \Tilde{Y}(k+T_p|k) = \Tilde{C} x(k) + \Tilde{D_u} U(k+T_c) + \Tilde{D_d} W(k+T_p),\label{eq:state_space_pred_2}
    \end{gather}
\normalsize
where the predictive system matrices $\Tilde{A},~\Tilde{B_u},~\Tilde{B_d}, ~\Tilde{C}$ and $\Tilde{D}$ can be computed using the definitions above.

In our formulation, we consider that the measured disturbance $w(k)$ in \eqref{eq:state_space_pred_1}-\eqref{eq:state_space_pred_2} remains constant within the prediction horizon $T_p$. Therefore, we have $ w(k+n|k) = w(k), \quad n=1,\ldots, T_p$.
Consequently, the disturbance vector can be computed as ${W}(k+T_p)=[w(k), \cdots, w(k)]^T$. The inaccuracy in modeling the predicted disturbance vector ${W}(k+T_p)$ can be compensated by incorporating a feedback scheme into the receding horizon optimization \cite{borrelli2017predictive}.

\subsection{Objective Functions}
Let us define $\left \| z \right \|_M$ to be the $M$ weighted norm of an arbitrary vector $z$ such that $\left \| z \right \|_M:=(z^TMz)^{\frac{1}{2}}$.
In order to drive each HDV's state towards the equilibrium platoon state, the primary aim of the CAV controller is to minimize the squared error between the predicted output $\Tilde{y}(k+n|k)$, $n=1,2, \ldots, T_p$, and the corresponding reference output. The first objective function thus takes the form $J_1 := \frac{1}{2} \sum_{n=1}^{T_p} \left \| \Tilde{y}(k+n|k))-y_r(k+n) \right \|^2_Q$,
where the reference output $y_r(k):=[0, ~(N-1)(s_0 + \rho \begin{bmatrix}
   1 &0
   \end{bmatrix}w(k)),~ 0]^T$ and the positive semi-definite output weight matrix $Q:=diag(q_v,q_{e_1},q_{e_2})$ 
with the diagonal weight parameters $q_v,~q_{e_1}, ~q_{e_2}$ corresponding to the speed $v_1(k)$, head-to-tail gap $e_{1,1}(k)$ and leader-follower gap $e_{1,2}(k)$, respectively. 
Since the measured disturbance $w(k)$ remains constant within the prediction horizon $T_p$, and the reference output $y_r(k)$ is an explicit function of the measured disturbance $w(k)$, the predictive reference output $y_r(k+n|k), ~n=1,\ldots,T_p$ remains constant within the prediction horizon $T_p$ as well. Thus we have
$y_r(k+n|k) = y_r(k), \quad n=1,\ldots, T_p.$

The second objective of the controller is to improve the fuel economy of the CAV by minimizing the $L^2$-norm of the CAV's control input. Hence, we have the second objective function $ J_2 := \frac{1}{2} \sum_{m=1}^{N_c}  \left \| u(k+m-1) \right \|^2_R$,
where $R:=[w_r]$ is the positive definite weight matrix on the control input with positive weight parameter $w_r$. 

Finally, combining the above objective functions and using the compact notations from \eqref{eq:state_space_pred_1}-\eqref{eq:state_space_pred_2}, we have the final objective function as follows
\small
\begin{align}
    J =  \frac{1}{2} \left \| \Tilde{Y}(k+T_p|k)-Y_r \right \|^2_{\bar{Q}}+\frac{1}{2}  \left \| U(k+T_c) \right \|^2_{\bar{R}},\label{eq:final_cost}
\end{align}
\normalsize
where $Y_r =[y_r(k), \cdots, y_r(k)]^T$, and $\bar{Q}$ and $\bar{R}$ are weight matrices.
\subsection{Constraints}
%
In our formulation, we consider the constraints on the control input in \eqref{eq:state_control_constraints}, safety in \eqref{eq:rearend_constraint}, and CAV speed in \eqref{eq:state_control_constraints} associated with the physical limitation of the CAV dynamics, passenger safety, and speed limit of the roadway, respectively.
The constraints in the context of the proposed receding horizon control framework are thus given as
\small
\begin{subequations}\label{eq:ocp_constraints}
\begin{align}
    &u_{\min} \le u(k+m-1) \le u_{\max},~ m=1,\ldots,T_c, \\
    &e_{1,1}(k+n) \ge (N-1)s_0, ~ n=1,\ldots,T_p,\\
    &e_{1,2}(k+n) \ge s_0, ~ n=1,\ldots,T_p,\\
    &v_{\min} \le v(k+n) \le v_{\max}, ~ n=1,\ldots,T_p.
\end{align}
\end{subequations}
\normalsize
\subsection{The Optimal Control Problem}
With the objective function \eqref{eq:final_cost}, constraints \eqref{eq:ocp_constraints}, dynamics model \eqref{eq:discrete_dynamics_1}-\eqref{eq:discrete_dynamics_4}, and the information set $\mathcal{I}_1(k),~k=0,\ldots, T_h$ at hand, the optimal control problem can finally be written as
\small
\begin{align}\label{eq:ocp}
&\min_{U(k+T_c)}  J,\\ 
&\text{subject to}: \eqref{eq:discrete_dynamics_1}-\eqref{eq:discrete_dynamics_4},\eqref{eq:ocp_constraints} \text{ and given } \mathcal{I}_1(k)\nonumber
\end{align}
\normalsize
The optimal control problem in \eqref{eq:ocp} can be transformed into a standard quadratic programming problem and solved using the active-set algorithm, see \cite{bemporad2004MPCMatlab,borrelli2017predictive}. It is possible to soften the state constraints in \eqref{eq:ocp_constraints} to facilitate the feasibility of the solution of \eqref{eq:ocp}. However, significantly large penalty should be incorporated into the objective function in \eqref{eq:ocp} using a dimensionless, non-negative slack variable to handle the soft constraint violation, the exposition of which is outside the scope of this paper and can be found in \cite{bemporad2004MPCMatlab}. 
\section{Simulation results}\label{sec:sim}
\begin{figure}
    \centering
    \includegraphics[scale=0.56]{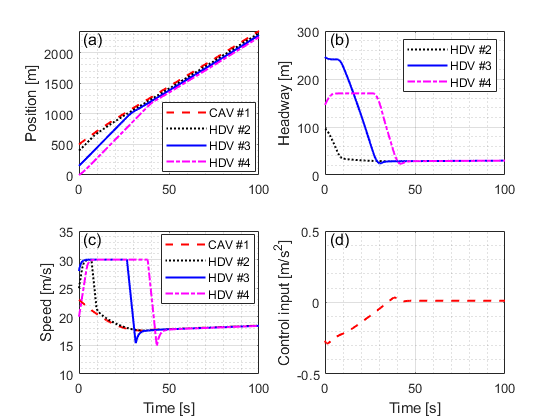}
    \caption{Platoon formation with OVM car-following model \eqref{eq:hdv_dynamics_ovm} for $N=4$, where the (a) position trajectory, (b) vehicle headway, (c) speed trajectory and (d) the CAV control trajectory are shown.}
    \label{fig:ovm_sim_1}
\end{figure}

To evaluate the performance of the proposed control framework, we adopt the optimal velocity model (OVM) \cite{bando1995dynamical} and the intelligent driver model (IDM) \cite{treiber2013traffic} to represent the predecessor-follower coupled dynamics of each HDV $i\in \mathcal{N}_{\text{HDV}}$. 
One of the simplest forms of the OVM car-following model \cite{bando1995dynamical} is given as
\small
\begin{align}
    &{{u}_i(t) = \alpha (V_i(\delta_i(t),s_i(t)) -v_i(t)),}\label{eq:hdv_dynamics_ovm} \quad \\
   &V_i(\delta_i(t),s_i(t))=
\begin{array}
[c]{ll}%
 {\frac{v_d}{2}(\tanh(\delta_i(t))}{+\tanh(s_i(t))),}
\end{array}\nonumber
\end{align}
\normalsize
where $\delta_i(t):= \Delta p_i(t) - s_i(t)$, and $\alpha$, $V_i(\delta_i(t),s_i(t))$ and $v_d$ denote the control gain representing the driver's sensitivity coefficient, the equilibrium speed-headway function 
and the desired speed of the roadway, respectively.
The IDM car-following model \cite{treiber2013traffic} for HDV $i\in\mathcal{N_{\text{HDV}}}$ has the following structure
\small
\begin{align}
    &u_i(t) = a\bigg( 1-\bigg(\frac{v_i(t)}{v_{d}}\bigg)^\gamma - \bigg(\frac{\Delta \bar{p}_i(t)}{\Delta p_i(t)}\bigg)^2 \bigg),\label{eq:hdv_dynamics_idm}\\
    &\Delta \bar{p}_i(t) = s_i(t) + \frac{v_i(t) \Delta v_i(t)}{2ab}\nonumber,
\end{align}
\normalsize
where, $a,~b$ and $\gamma$ are the desired acceleration, comfortable braking and acceleration exponent, respectively.
The parameters for the car-following models and the receding horizon controller considered in our numerical study can be found in \href{https://sites.google.com/view/ud-ids-lab/model-agnostic-platoon}{https://sites.google.com/view/ud-ids-lab/model-agnostic-platoon}.
We conduct the simulation studies using MATLAB R2020b/Simulink with the configuration of Intel Core i7-6700 CPU @ 3.40 GHz.
%
%
For the first case study, a platoon formation for $N=4$ vehicles is shown Fig. \ref{fig:ovm_sim_1}, where the OVM model in \eqref{eq:hdv_dynamics_ovm} is considered for the trailing HDVs. The leading CAV and trailing HDVs have randomly selected initial position (Fig. \ref{fig:ovm_sim_1}(a)) and initial speed (Fig. \ref{fig:ovm_sim_1}(c)), respectively. The lead CAV implements the proposed controller to complete the platoon formation operation near $50$ s (according to Remark \ref{rem:rmse_platoonFormation}), and the vehicle headway (Fig. \ref{fig:ovm_sim_1}(b)) and speed (Fig. \ref{fig:ovm_sim_1}(c)) converge to some equilibrium value. Additionally, none of the constraints in \eqref{eq:ocp_constraints} were violated as evident from the headway profile in Fig. \ref{fig:ovm_sim_1}(b), speed profile in Fig. \ref{fig:ovm_sim_1}(c), and CAV's control input trajectory in Fig. \ref{fig:ovm_sim_1}(d), respectively.
\begin{figure}
    \centering
    \includegraphics[scale=0.5]{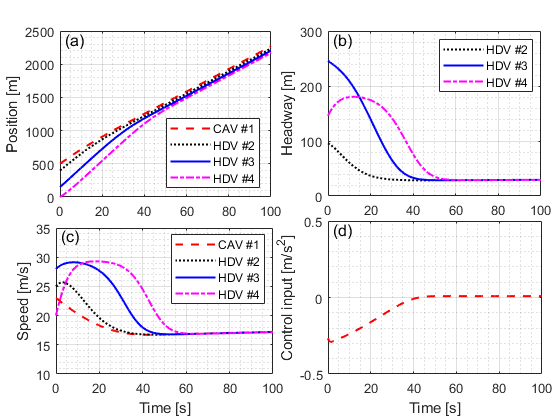}
    \caption{Platoon formation with IDM car-following model \eqref{eq:hdv_dynamics_idm} for $N=4$, where the (a) position trajectory and (b) vehicle headway, (c) speed trajectory and (d) the CAV control trajectory are shown.}
    \label{fig:idm_sim_1}
\end{figure}
To validate the model-agnostic nature of the proposed controller, we present a second case study using the IDM model \eqref{eq:hdv_dynamics_idm} (see Fig. \ref{fig:idm_sim_1}) considering the same initial conditions as in the previous case, which yields similar result without violating any constraints in \eqref{eq:ocp_constraints}, as shown in Figs. \ref{fig:idm_sim_1}(a)-(d). It is interesting to note that, we have monotonically increasing non-positive linear optimal control input trajectory of the CAV in both of the above cases (see Figs. \ref{fig:ovm_sim_1}(d) and \ref{fig:idm_sim_1}(d)), which resembles a typical energy-optimal control input trajectory derived using standard Hamiltonian analysis \cite{mahbub2020Automatica-2, Malikopoulos2020}. Note that, we can consider a mixture of OVM and IDM car-following model for the HDVs by appropriately selecting $\epsilon_{\Delta p}$ and $\epsilon_{v}$ in \eqref{eq:platoon_formation_rmse}.
\begin{figure}
    \centering
    \includegraphics[scale=0.55]{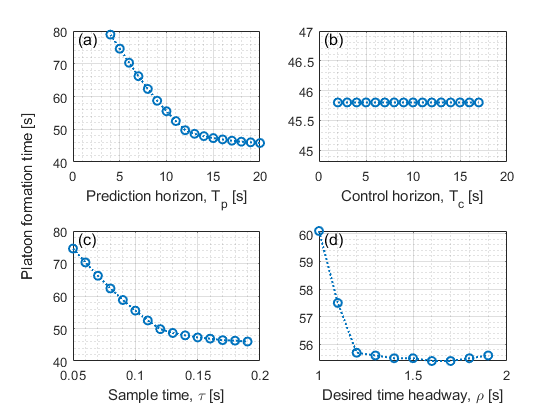}
    \caption{Platoon formation with IDM car-following model \eqref{eq:hdv_dynamics_idm} for $N=4$, where the sensitivity of the platoon formation time under varying (a) prediction horizon $T_p$, (b) control horizon $T_c$, (c) sample time $\tau$ and (d) desired time headway $\rho$ are illustrated.} 
    \label{fig:sensitivity_1}
\end{figure}
\begin{figure}
    \centering
    \includegraphics[width=3.4in]{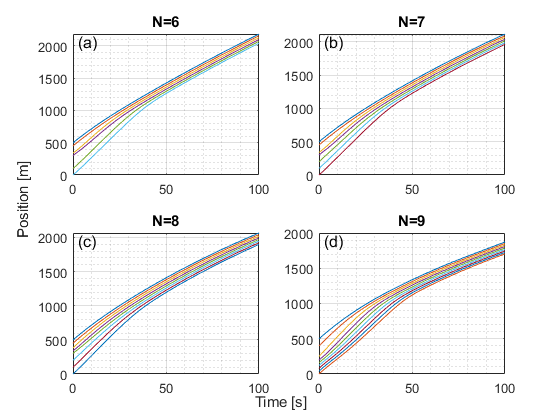}
    \caption{Safe platoon formation with IDM car-following model \eqref{eq:hdv_dynamics_idm} for $6,7,8$ and $9$ vehicles.} 
    \label{fig:robustness_1}
\end{figure}

Figure \ref{fig:sensitivity_1} shows the sensitivity analysis of the proposed control framework for $N=4$ subject to varying controller parameters $T_p, T_c, \text{ and } \tau$, and IDM car-following parameter $\rho$. Here, we use \eqref{eq:platoon_formation_rmse} to compute the platoon formation time. Increasing $T_p$ and $\tau$ decrease the platoon formation time, as shown in Figs. \ref{fig:sensitivity_1}(a) and \ref{fig:sensitivity_1}(c), respectively, whereas the variation of $T_c$ does not affect the platoon formation time, as shown in Fig. \ref{fig:sensitivity_1}(b). However, choosing appropriate $T_c$ is essential to enforce the constraints in \eqref{eq:ocp_constraints}. Note that, the parameters $T_p$ and $\tau$ can be tuned using Figs. \ref{fig:sensitivity_1}(a) and \ref{fig:sensitivity_1}(c) to form a platoon within the desired optimization horizon.
The platoon formation time under varying $\rho$, which represents different driving behavior of the IDM model, is shown in Fig. \ref{fig:sensitivity_1}(d). Here, the proposed framework is robust enough to form a platoon within the optimization horizon $T_h=65$ s. In all of the cases presented in Figs. \ref{fig:sensitivity_1}(a)-(d), the proposed controller enables platoon formation without violating any constraints in \eqref{eq:ocp_constraints}. 
Finally, we investigate the robustness of the proposed framework under different platoon size $N=6,~7,~8$ and $9$ as shown in Fig. \ref{fig:robustness_1}. The position trajectories in Fig. \ref{fig:robustness_1} indicates that the CAV controller is able to form platoon within the optimization horizon $T_h=65$ s without violating any safety constraint in \eqref{eq:rearend_constraint}.
\section{Discussion and concluding Remarks}\label{sec:conc}
In this paper, we presented a constrained multi-objective optimal control framework for platoon formation under a mixed traffic environment, where a leading CAV computes and implements its optimal control input to force the following HDVs to form a platoon. We developed a model-agnostic receding horizon control framework with a multi-successor communication topology that solves in real time the optimal control problem, and provided detailed sensitivity and robustness analysis using numerical simulation to validate the performance of the proposed framework.

{A direction for future research should extend the proposed framework for optimal coordination of mixed vehicle platoon in traffic scenarios such as on-ramp merging, urban intersection, etc. Ongoing research investigates the incorporation of non-linear state-space representation and different communication topology to improve the controller performance.}
\bibliographystyle{IEEEtran}
\bibliography{references/IDS_Publications_09082021,references/acc_pt_vd_ref, references/platoon, references/mpc,references/TAC_references}

%

%

\end{document}